\documentclass[11pt]{article}

\input{amssym.def}
\input{amssym}
\usepackage{eufrak}

\newtheorem{theorem}{Theorem}

\newtheorem{conjecture}[theorem]{Conjecture}

\newtheorem{definition}[theorem]{Definition}

\newtheorem{proposition}[theorem]{Proposition}
\newtheorem{remark}[theorem]{Remark}

\def\beq{\begin{equation}}
\def\be{\begin{equation}}
\def\eeq{\end{equation}}
\def\ee{\end{equation}}
\def\K{{\Bbb K}}
\def\C{{\Bbb C}}
\def\R{{\Bbb R}}
\def\Z{{\Bbb Z}}

\def\gg{{\frak g}}

\def\PP{{\cal P}}

\def\NN{{\Bbb N}}
\def\D{{\frak D}}

\def\N{{\cal N}}
\def\A{{\cal A}}
\def\RR{{\cal R}}

\def\D{{\cal D}}
\def\W{{\cal W}}

\def\q{q^{-1}}

\def\pa{\partial}
\def\dd{\tilde{\partial}_t}

\def\tpa{\tilde{\pa_t}}
\def\V{\cal V}

\def\LR{{\cal L}(R)}
\def\gh{{\frak g}_{\h}}

\def\Uq{U_q(sl(m))}

\def\gg{\mbox{$\frak g$}}
\def\ot{\otimes}
\def\h{\hbar}

\def\cas{{\rm cas}}

\def\Ob{{\rm{Ob}}}
\def\Mor{{\rm{Mor}}}
\def\vl{V_\la}

\def\vvm{V^*_\mu}

\def\vv{V^{\ot 2}}

\def\Sym{{\rm Sym}}
\def\Im{{\rm Im}}
\def\Tr{{\rm Tr}}

\def\End{{\rm End}}

\def\det{{\rm det}}

\def\diag{{\rm diag}}

\def\De{{\Delta}}
\def\de{{\delta}}

\def\Cas{{\rm Cas}}

\def\al{\alpha}

\def\la{{\lambda}}

\def\A{{\cal A}}
\def\Ah{{\cal A}_\h}

\def\hh{\displaystyle\frac{\h}{2}}

\begin{document}

\makeatletter
\renewcommand{\theequation}{{\thesection}.{\arabic{equation}}}
\@addtoreset{equation}{section}
\makeatother

\title{Braided algebras and their applications to Noncommutative Geometry}
\author{
\rule{0pt}{7mm} Dimitri
Gurevich\thanks{gurevich@univ-valenciennes.fr}\\
{\small\it LAMAV, Universit\'e de Valenciennes,
59313 Valenciennes, France}\\
\rule{0pt}{7mm} Pavel Saponov\thanks{Pavel.Saponov@ihep.ru}\\
{\small\it
National Research University Higher School of Economics,}\\
{\small\it 20 Myasnitskaya Ulitsa, Moscow 101000, Russia}\\
{\footnotesize\it \&}\\
{\small\it Division of Theoretical Physics, IHEP, 142281
Protvino, Russia} }

\maketitle

\begin{abstract}
We introduce the notion of a braided algebra and study some examples of these. In particular, $R$-symmetric and $R$-skew-symmetric algebras of a linear space $V$ equipped with a skew-invertible Hecke symmetry $R$ are braided
algebras. We prove the ``mountain property" for the numerators and denominators of their Poincar\'e-Hilbert series (which
are always rational functions).

Also, we further develop a  differential calculus on  modified Reflection Equation algebras. Thus, we exhibit a new form of
the Leibniz rule for partial derivatives on such algebras related to involutive symmetries $R$. In particular, we present this
rule for the algebra $U(gl(m))$. The case of the algebra $U(gl(2))$ and its compact form $U(u(2))$ (which can be treated
as a deformation of the Minkowski space algebra) is considered in detail. On the algebra $U(u(2))$ we introduce the notion
of the quantum radius, which is a deformation of the usual radius, and compute the action of rotationally invariant operators
and in particular of the Laplace operator. This enables us to define analogs of the Laplace-Beltrami operators corresponding
to certain Schwarzschild-type metrics and to compute their actions on the algebra $U(u(2))$ and its central extension. Some
"physical" consequences of our considerations are presented.
\end{abstract}

{\bf AMS Mathematics Subject Classification, 2010:} 17B37, 81R60


\section{Introduction}

The aim of this paper is twofold. On the one hand, we introduce the notion of a braided algebra and further develop analysis
on the Reflection Equation (RE) algebra, which is an example of such an algebra. On the other hand, we exhibit certain
applications of the RE algebras to $U(u(2))$-covariant Noncommutative (NC) Geometry.

By braided algebras we mean a subclass of unital algebras related to braidings, whereas the term braiding stands for an
invertible operator $R:\vv\to\vv $ ($V$ is a given finite dimensional linear space) satisfying the so-called Quantum Yang-Baxter Equation
$$
R_{12}\,R_{23}\,R_{12}=\,R_{23}\,R_{12}\,R_{23},\qquad R_{12}=R\ot I,\quad R_{23}= I\ot R.
$$
Hereafter, $I\in \End(V)$ is the identity operator and the lower indices indicate the factors of the tensor product on which
a given operator acts.

Below, we are dealing with braidings satisfying an additional condition
$$
(R-q I)\, (R+\q I)=0, \quad q\in\K,
$$
where $\K$ is the ground field ($\C$ or $\R $) and the parameter $q$ is assumed to be generic. Such braidings are called
{\it Hecke symmetries}. For example, Hecke symmetries come from the quantum groups $\Uq $. These Hecke symmetries
(and other objects related to $\Uq $) will be referred to as {\it standard} ones. Note, that quantum groups of the $B$, $C$
and $D$ series provide braidings of the Birman-Murakami-Wenzl type.

Our definition of braided algebras is given in terms of objects and morphisms of the so-called Schur-Weyl category $SW(V)$
generated by a given linear space $V$ equipped with a skew-invertible Hecke symmetry $R:\vv\to\vv $ (see section
\ref{sec:2}). Typical examples of braided algebras are the $R$-symmetric, and the $R$-skew-symmetric algebras of the
space $V$ and the Reflection Equation (RE) algebra associated with a given Hecke symmetry $R$. By contrast, the
well-known RTT algebra (which in the standard case is the {\it restricted dual} Hopf algebra to $\Uq $) is not a braided
algebra.

All aforementioned braided algebras have a good deformation property. For quadratic algebras this means that for a generic
$q$ the dimensions of their homogeneous components equal those for $q=1$. As for the RE algebra, it has another important
property: it admits a change of generators converting it into a quadratic-linear algebra ({\it modified} RE algebra) similar to an enveloping algebra. More precisely, in the standard case the modified RE algebra tends to the enveloping algebra $U(gl(m))$
as $q\to 1$. Also, if a Hecke symmetry $R$ is a deformation of the super-flip on the super-space $V_{m,n}$, the corresponding modified RE algebra tends to the super-algebra $U(gl(m|n))$ as $q\to 1$. Besides, the RE algebra corresponding to a proper Hecke symmetry has a representation category looking like that of the algebra $U(gl(m))$ (or $U(gl(m|n))$) and turns into the
latter one as $q\to 1$ (see \cite{GPS1}). Due to this fact, generators of the RE algebra are good candidates for the role of
vector fields. Such vector fields were introduced in the framework of quantum differential calculus on a pseudogroup (in fact, on the RTT algebra or on its compact version) in the series of papers \cite{W, IP, FP}.

In \cite{GPS5} we suggested a more general calculus in which the role of the function algebra Fun$(GL(m))$ is played by a quantum matrix algebra while the role of vector fields is always played by a (modified) RE algebra. The most interesting
particular case arises when the role of the algebra Fun$(GL(m))$ is taken by another copy of the RE algebra. In this case we
get a braided analog of the Heisenberg double\footnote{The conventional Heisenberg double consists of two Hopf algebras dual to each other.} constructed from two copies of the RE algebra. In the paper \cite{GPS6} we used such a "braided Heisenberg double" in order to define  braided analogs of partial derivatives on the RE algebra playing the role of Fun$(GL(m))$. Besides,
by passing to the limit $q\to 1$ in the standard case, we constructed the partial derivatives on the NC algebras $U(gl(m))$.

These derivatives differ from the usual ones, defined on the commutative algebra $\Sym(gl(m))$, by a modification of the
Leibniz rule. Below, we exhibit a form of this rule well adapted to algebras $U(gl(m))$ and their generalized analogs associated
with any skew-invertible involutive symmetry\footnote{\label{ftn:inv}We say that a braiding $R$ is an involutive symmetry if $R^2=I$.} $R$. Our construction  includes the differential calculus (with the standard Leibniz rule) on the super-algebras $U(gl(m|n))$. Thus, by making use of braided algebras we get some elements of NC Geometry on algebras having no "braided
ingredient" in their construction, the latter algebras are covariant with respect to the usual groups (or their super-analogs).

In the particular case $m=2$, $n=0$ we get a version of the differential calculus on the algebra $U(u(2))$\footnote{In fact, we deal with the algebra $U(u(2)_\h)$ where the notation $\gh $ means that a deformation multiplier $\h $ is introduced in the Lie bracket of a given Lie algebra $\gg $.}. An amazing fact is that  partial derivatives defined on the algebra $U(u(2)_\h)$ commute with each other. Consequently, any classical differential operator with constant coefficients can be easily re-defined on this
algebra. Thus, analogs of the Klein-Gordon, Maxwell and Dirac operators on the algebra $U(u(2)_\h)$\footnote{This algebra can be seen as a NC deformation of the Minkowski space algebra, though it does not admit any reasonable analog of the Lorentz  subgroup.} have forms just the same as the classical ones but with a new meaning of the partial derivatives.

The situation is more complicated for operators with coefficients from $\Sym(u(2))$. In order to transfer such an operator to
the algebra $U(u(2)_\h)$ we have to use a "quantizing map" $\Sym(u(2))\to U(u(2)_\h)$.  Such a map enables us to quantize
a differential operator with coefficients from $\Sym(u(2))$ to an operator with coefficients from $U(u(2)_\h)$ and acting on this algebra.

Nevertheless, many interesting (and physically meaningful) operators have coefficients belonging to a larger algebra. As an example of such an algebra $\A $ we consider the product $\K(t, r)\ot \Sym(u(2))$ where $\K(t, r)$ is the algebra of rational functions in the time $t$ and the radius $r=\sqrt{x^2+y^2+z^2}$ (hereafter $x,\,y,\,z$ are spatial variables). In particular, the Laplace-Beltrami operator corresponding to the Schwarzschild metric is of this form. We introduce a quantum counterpart $\Ah$
of the algebra $\A$ by defining a quantum analog $\hat r$ of the radius $r$. It should be emphasized that the "quantum radius" $\hat r$ is naturally extracted from the Cayley-Hamilton identity valid for the generating matrix of the algebra $U(u(2)_\h)$. (In
a more general context, such an identity is valid for the generating matrices of the RE algebras, see \cite{GPS2, GPS3,GPS4}.) Then, we extend the quantizing map to a map $\al:\A\to \Ah$.

Let $\D $ be a differential operator with coefficients from $\A $. By applying the quantizing map $\al $ to its coefficients we get the differential operator denoted $\al(\D)$. Its coefficients belong to the algebra $\Ah $ and the derivatives $\pa_t,\,...,\,\pa_z$ appearing in this operator are assumed to act on the algebra $U(u(2)_\h)$. We show how to extend some operators of this
form to the algebra $\Ah$: in particular we do this for the Laplace-Beltrami operator, corresponding to the Schwarzschild metric.

Note that the quantized operators we are dealing with, being expressed via the commutative variables $t$ and $\hat r$,
become difference operators. Thus, in contrast with usual differential operators, we can consider them on a lattice, i.e. on a discrete space-time. In a sense, its discreteness  is a consequence of its noncommutativity. Another ``physical" observation is that the space and time on the algebra $U(u(2)_\h)$ do not exist independently of one another. Indeed, the partial derivative
in $t$ involves an action on the spatial variables and vice versa. It seems to be plausible that a similar NC calculus does not
exist on the algebra $U(su(2)_\h)$, i.e. on the NC space without the time.

The paper is organized as follows. In the next section we introduce the notion of a braided algebra and consider some basic examples. In particular, we go back to the question of the possible form of a Hecke symmetry and prove the {\it mountain property} for the numerator and denominator of the Poincar\'e-Hilbert series corresponding to the (skew-)symmetric algebra defined via a skew-invertible Hecke symmetry. In section 3 we introduce Weyl algebras defined on RE algebras and consider their $q\to 1$ limits. Here, our main objective is to give a version of the Leibniz rule useful for further applications. In section 4
we introduce Laplace and other rotationally invariant operators acting on the RE algebras in question. We treat the case of the algebra $U(u(2)_\h)$ (which is a limit of a standard modified RE algebra) in detail. On this basis, in section 5 we construct a NC analog of a model describing the dynamics of a scalar massless field in a space equipped with a Schwarzschild-type metric.
We consider this example as the starting point for a study of higher dimensional (matrix) models. Also note, that our approach
is valid for the RE algebras (modified or not) but computations become much more complicated.
\bigskip

\noindent
{\bf Acknowledgement}
This study was carried out within The National Research University Higher School of Economics Academic Fund Program in 2012-2013, research grants No. 11-01-0042 and No. 12-09-0064. Also, the work of P.S. was partially supported by the RFBR grant 11-01-00980-a.

\section{Braided quadratic algebras: basic examples}
\label{sec:2}

Let $R=R(q):\vv\to \vv $ be a Hecke symmetry analytically depending on the parameter $q$ in a neighborhood $\V $
(assumed to be connected) of 1. Note that the symmetry $R(1)$ is involutive. Let us associate with this Hecke symmetry two quotient algebras of the free tensor algebra $T(V)$ of the space $V$: $R$-symmetric $\Sym_R(V)$ and $R$-skew-symmetric algebra $\bigwedge_R(V)$ defined as follows
\be
\Sym_R(V)=T(V)/\langle \Im (q\,I-R)\rangle,\qquad {\bigwedge}_R(V)=T(V)/\langle \Im (\q\,I+R)\rangle.
\label{sym-wedge}
\ee
Hereafter, the symbol $\langle X\rangle $ stands for the two-sided ideal generated by a subset $X$.

Though the complete classification of Hecke symmetries is an open problem, some information concerning their possible forms
can be drawn from the Poincar\'e-Hilbert (PH) series of the algebras introduced above
$$
\PP_+(t)=\sum_{k} \dim \Sym^{(k)}_R(V)\, t^k,\qquad
\PP_-(t)=\sum_{k} \dim {\bigwedge}^{(k)}_R(V)\, t^k,
$$
where $A^{(k)}$ stands for the degree $k$ homogenous component of a graded algebra $A$.

Constructed in the paper \cite{G2}, were two series of projectors $P_{\pm}^{(k)}$ which are $q$-analogs of the usual
symmetrization and skew-symmetrization operators:
$$
P_+^{(k)}:\,\,V^{\ot k}\to \Sym^{(k)}_R(V),\qquad P_-^{(k)}:\,\,V^{\ot k}\to {\bigwedge}^{(k)}_R(V),
\quad k\ge 2.
$$
Namely, if the parameter $q$ is subject to the condition
\be
n_q= {{q^n-q^{-n}}\over{q-\q}}\not=0,\quad n\in {\Bbb N}
\label{cond}
\ee
then there are isomorphisms of vector spaces
$$
\Sym^{(k)}_R(V)\cong \Im\, P_+^{(k)},\qquad
{\bigwedge}^{(k)}_R(V)\cong \Im\, P_-^{(k)}.
$$

Since $\dim\, \Sym^{(k)}_R(V)$ and $\dim\, {\bigwedge}^{(k)}_R(V)$ are lower semi-continuous functions in $q$ while the
functions $\dim \, \Im\, P_\pm^{(k)}$ are upper semi-continuous, we conclude that these functions are constant on the
connected neighborhood $\V$ with the possible exception of the values of $q$ violating the condition (\ref{cond}). Moreover,
as was shown in \cite{G2} for a generic $q$ (i.e. for all $q$ excepting a countable set not including 1) the relation
$$
\PP_+(t) \, \PP_-(-t)=1
$$
is valid. Hereafter, we consider the PH series $\PP_\pm(t)$ at generic values of the parameter $q$.

Let us recall, that a braiding $R$ is called skew-invertible if there exists an endomorphism
$\Psi: V^{\otimes 2}\rightarrow V^{\otimes 2}$ such that
\be
\Tr_{(2)} \,R_{12}\,\Psi_{23}=P_{13}=\Tr_{(2)} \, \Psi_{12}\,R_{23}.
\label{psi}
\ee
Hereafter, $P$ stands for the usual flip. Note that if a Hecke symmetry $R(q)$ is skew-invertible for a value of the parameter
$q$, it is skew-invertible for all generic $q$.

A skew-invertible braiding $R:\vv\to\vv $ can be extended to a braiding
\be
R:\,\,(V\oplus V^*)^{\ot 2} \to (V\oplus V^*)^{\ot 2},
\label{ex}
\ee
where we keep for the extended braiding the same notation $R$. In the above formula the space $V^*$ is dual to $V$, i.e. the
spaces $V$ and $V^*$ are equipped with a nondegenerated $R$-invariant pairing $\langle\,,\,\rangle:\,\, V\ot V^*\to\K $. The
pairing is $R$-invariant if the following property holds true (below the space $W$ stands for either $V$ or $V^*$)
\be
\begin{array}{rcr}
R \langle\,,\,\rangle_{23}=\langle\,,\,\rangle_{12} R_{23} R_{12}&\quad {\rm on}\quad &W\otimes V\otimes V^*\\
\rule{0pt}{5mm}
R \langle\,,\,\rangle_{12}=\langle\,,\,\rangle_{23} R_{12} R_{23}&\quad {\rm on}\quad & V\otimes V^*\otimes W,
\end{array}
\label{pair}
\ee
where we assume that the ground field $\K $ commutes with $V$ and $V^*$ in the usual way:
$$
R(w\ot k)=k\ot w,\qquad R(k\ot w)=w\ot k,\qquad \forall\, k\in\K,\,\, w\in W.
$$

The prolongation (\ref{ex}) exists and is unique. In order to describe it explicitly, we fix dual basis $\{x_i\}$ and $\{x^j\}$ in the
spaces $V$ and $V^*$ respectively:
$$
\langle x_i, x^j \rangle= \de_i^j.
$$

\begin{remark} {\rm Note that in the space $V^*$ one can define right and left dual basis to a given basis of $V$. We are
dealing with the right dual basis.}
\end{remark}

Let $\{x_i\ot x_j\}$ be the corresponding basis of the space $\vv $. In this basis the skew-invertible Hecke symmetry
$R: \vv\to\vv $ is represented by a $(\dim V)^2 \times(\dim V)^2 $ matrix $\|R_{ij}^{kl}\|$
$$
R(x_i\ot x_j)=x_k\ot x_l\,R_{ij}^{kl}.
$$
Hereafter, summation over repeated indices is always assumed.

Then the extension (\ref{ex}) can be represented by the following matrices on the components $\vv $, $(V^*)^{\ot 2}$,
$V\ot V^*$ and $V^*\ot V$ of the tensor square $(V\oplus V^*)^{\otimes 2}$
\be
\begin{array}{lc}
R(x_i\ot x_j)=x_k\ot x_l\,R_{ij}^{kl}, & R( x^i\ot x^j)=x^k\ot x^l\,R^{ji}_{lk},\\
\rule{0pt}{6mm}
R(x_i\ot x^j)=x^k\ot x_l\,(R^{-1})_{ki}^{lj}, & R(x^j\ot x_i)= x_k\ot x^l\,\Psi_{li}^{kj}.
\end{array}
\label{ext}
\ee

It is easy to show that the embedding $\K\to V^*\ot V$ generated by the correspondence
\be
1\mapsto \sum_i x^i \ot x_i
\label{emb}
\ee
is $R$-invariant.

Let us introduce two endomorphisms $B$ and $C$ of the space $V$ via the morphism $\Psi$ (\ref{psi}) as $B:=\Tr_{(1)}\,\Psi$, $C:=\Tr_{(2)}\,\Psi\,$ or, in coordinate form
\be
B_i^j=\Psi_{ki}^{kj},\qquad C_i^j=\Psi_{ik}^{jk}. \label{BC}
\ee
Note that if $R$ is a super-flip, the operators $B$ and $C$ are equal to each other and coincide with the parity operator.

\begin{proposition} {\bf (\cite{H,Da, DH})}.
\label{prop2}
The HP series $\PP_-(t)$ (and hence $\PP_+(t)$) is a rational function:
{\rm
\be \PP_-(t)=\frac{N(t)}{D(t)}=
\frac{1+a_1\,t+...+a_m\,t^m}{1-b_1\,t+...+(-1)^n\, b_n\,t^n} =
\frac{\prod_{i=1}^m(1+x_it)}{\prod_{j=1}^n(1-y_jt)}\,,
\label{p-}
\ee}
$\!\!$where $a_i$ and $b_i$ are positive integers, the polynomials $N(t)$ and $D(t)$ are coprime, and all the numbers
$x_i$ and $y_i$ are real positive.

If, in addition, the Hecke symmetry is skew-invertible, then the polynomials $N(t)$ and $D(-t)$ are reciprocal.
\end{proposition}

Recall, that a polynomial $p(t) = c_0 + c_1t +\dots +c_nt^n$ is called {\it reciprocal} if $p(t) = t^np(t^{-1})$ or, equivalently,
$c_i = c_{n-i}$, $0\le i\le n$. Note, that if a real number $z\not=1$ is a root of a reciprocal polynomial, the number $z^{-1}$
is a root as well with the same multiplicity. This entails that the polynomials $N(t)$ and $D(-t)$ factorize into a product of terms
of the form
\be
(1+t)\quad {\rm and}\quad (1+z t)(1+z^{-1}t)=1+ct+t^2\quad {\rm where}\quad c\geq 2.
\label{prodd}
\ee
For such a polynomial the following claim is valid.
\begin{proposition} \label{prop:mount}
Any polynomial which is a product of terms listed in (\ref{prodd}) has the "mountain property": its coefficients strictly increase
up to the middle (and, consequently, strictly decrease after it).
\end{proposition}
The proof of the above proposition can be easily obtained by induction on the number of the factors in such a polynomial.

Proposition \ref{prop:mount} provides new information on the possible form of a Hecke symmetry. In particular, if a skew-invertible Hecke symmetry is even (see definition \ref{def:bi-rank} below) then its HP series $\PP_-(t)$ is a polynomial possessing the mountain property.

However, we do not know whether the coefficients $c$ in the multipliers (\ref{prodd}) entering the factorized form of $N(t)$
and $D(-t)$ must be integer: we do not know, for example, whether there exists a Hecke symmetry such that its HP series $\PP_-(t)$ is equal to
$$
1+10t+12t^2+10 t^3+t^4=(1+(5+\sqrt{15})t+t^2)(1+(5-\sqrt{15})t+t^2).
$$

Nevertheless, if the monic polynomials $N(t)$ and $D(-t)$ are products of factors (\ref{prodd}) with {\it integer} middle coefficients $c$, then there exists a Hecke symmetry $R(q)$ for which $\PP_-(t)=\frac{N(t)}{D(t)}$. It can be constructed
by methods of \cite{G2}.

\begin{definition}
\label{def:bi-rank}
{\rm
Let $R:\vv \to\vv$ be a skew-invertible Hecke symmetry. Let $m$ (respectively $n$) be the degree of the numerator
(respectively denominator) of the corresponding PH series $\PP_-(t)$. The couple $(m|n)$ is called the {\it bi-rank} of
the space $V$ or of the Hecke symmetry $R$. If the bi-rank is $(m|0)$ (respectively $(0|n)$) the corresponding Hecke
symmetry is called {\it even} (respectively odd).}
\end{definition}

This notion is a generalization of the super-dimension for super-spaces and equals the latter provided $R$ is a super-flip.

Now, consider the monoidal quasitensor rigid category introduced in the paper \cite{GPS1}. This category is generated by
the spaces $V$ and $V^*$. Following \cite{GPS1} we call it the {\it Schur-Weyl} category and denote by the symbol
$SW(V)$. We briefly describe  the category $SW(V)$ with minor modifications comparing with \cite{GPS1}.

Let
$$
\la=(\la_1, \la_2,...,\la_k),\qquad \la_1\geq \la_2\geq...\geq\la_k, \qquad \sum\, \la_i=k
$$
be a partition of an integer $k\geq 2$. There is known a Schur functor $V\to \vl $ which assigns a space
$\vl\subset V^{\ot k}$ to the {\it basic} space $V$. In \cite{GPS1} ``braided analogs" of the spaces $\vl $ were introduced
as images of some projectors $P_\la: V^{\ot k}\to \vl $ constructed via a given Hecke symmetry $R$. In general, there is a
family of equivalent projectors (and, consequently, isomorphic spaces $\vl $) which differ from each other by embeddings into
the space $V^{\ot k}$. By $\vl $ we mean any of these spaces. In a similar manner the subspaces $\vvm \subset (V^*)^{\ot k}$
can be introduced.

By definition, the class $\Ob(SW(V))$ of objects of the category $SW(V)$ includes all spaces $V_\lambda$,
$V^*_\mu$ for all possible partitions $\lambda$ and $\mu$ as well as all their tensor products and direct sums
of these products\footnote{Note that for any given Hecke symmetry $R$ with a bi-rank $(m|n)$ some of the
projectors $P_\lambda $ are zero operators (and, consequently, spaces $V_\lambda$ and $V^*_\lambda$ vanish).
If the Hecke symmetry $R$ is even or odd, the category $SW(V)$ can be introduced with the use of the only
space $V$ (see \cite{G2} for detail).}.

Now, we want to describe the family $\Mor(SW(V))$ of morphisms of this category. For this purpose, we
note that any object of the category can be given a comodule structure over the so-called RTT algebra.
This is an associative unital algebra generated by the elements $t_i^j$, $1\le i,j\le\dim V$, which are subject
to the quadratic relations
$$
R\,T_1\, T_2=T_1\, T_2\, R,\qquad T_1=T\ot I,\quad T_2=I\ot T, \quad T=\|t_i^j\|.
$$
For arbitrary $R$ this is a bialgebra, but for the standard $R$ the RTT algebra is a Hopf algebra\footnote{To
introduce the Hopf structure one needs an antipodal map in the bialgebra. For this purpose, the RTT bialgebra
should be extended by the element $(\det_{\!R} T)^{-1}$, where $\det_{\!R} T$ is the so-called quantum determinant.}
restricted dual to the quantum universal enveloping algebra $U_q(gl(m))$ $m=\dim V$ (see \cite{FRT} for detail).

In general, if $R$ is an even Hecke symmetry, this bialgebra structure can be also extended up to the Hopf one (see
\cite{GPS1} for detail). Then the comodule structure on the basis space $V$, on its dual $V^*$ and on
$V\otimes V^*\cong \End (V)$ reads
\be
x_i\mapsto t_i^k\ot x_k,\quad x^j\mapsto s_k^j\ot x^k,\quad
x_i\,x^j \mapsto t_i^ks_p^j\ot x_k\, x^p,
\label{coac}
\ee
where $s_k^j$ are entries of the "inverse matrix" $S(T)$. Here $S$ is the antipode. The comodule structure on an arbitrary
$U\in \Ob(SW(V))$ is now easily defined since the above coaction can be naturally extended up to any tensor products of $V$ and $V^*$ and it commutes with the projectors $P_\lambda$.

By definition, a linear map $\beta : U_1\to U_2$, $U_1,U_2\in \Ob(SW(V))$ belongs to $\Mor(SW(V))$ iff it
commutes with the coaction of the RTT algebra on the spaces $U_1$ and $U_2$. In the standard case the
RTT coaction can be changed for the action of the quantum group $U_q(gl(m))$. So, in this case a morphism is,
by definition, a linear map commuting with the action of the quantum group.

For instance, it turns out that the initial skew-invertible braiding $R$ can be naturally extended up to a family of linear maps
$$
R_{U_1,U_2}:\,\, U_1\ot U_2\to U_2\ot U_1,\quad \forall\,U_1,\, U_2\in \Ob(SW(V)).
$$
All these maps commute with the RTT coaction and form the family of {\it permutation morphisms} of the
Schur-Weyl category.

There exists another way of introducing the set $\Mor(SW(V))$ (see \cite{GPS1}) but that given above is technically
more useful.

\begin{definition}
\label{def:6}{\rm
Algebras of the form
$$
T(U)/\langle \beta(W)\rangle,\qquad U,\,W\in \Ob(SW(V)),
$$
where $W$ is a subspace of the free tensor algebra $T(U)$ and $\beta:W\to W$ is an element of $\Mor(SW(V))$
are called braided. The relations of the form $\beta(W)=0 $ are called {\it admissible}.
}
\end{definition}

It is clear that the algebra $\Sym_R(V)$ (respectively $\bigwedge_R(V)$) is braided. To show this we set
$U=V$, $W=\vv$ and $\beta=q\, I-R$ (respectively $\beta=\q\, I+R$).

\begin{remark}
\label{rem7}{\rm
If $U\in \Ob(SW(V))$ is an object, then the direct sum $U\oplus U$ is an object as well. In order to distinguish these two
copies of the object $U$ we denote one of them by $U'$. Let $W\subset T(U\oplus U')$ be a subspace of the form
$(U\ot U)\oplus(U'\ot U')\oplus (U\ot U'\oplus U'\ot U)$. Let $\beta $ be a map such that its restrictions to each of the three
above components are morphisms. Then it is a morphism and the corresponding braided algebra is defined by three sets of
relations:
$$
\beta(U\ot U)=0,\quad \beta(U'\ot U') = 0, \quad \beta(U\ot U'\oplus U'\ot U)=0.
$$
Below, we consider such algebras constructed from two copies of an object and a morphism $\beta $ of the indicated
form. The third set of relations plays the role of permutation relations between the algebras
$T(U)/\langle \beta(U\ot U)\rangle$ and $T(U')/\langle \beta(U'\ot U')\rangle$.
}
\end{remark}

Now, we consider the object $\End(V) \cong V\ot V^*\in \Ob(SW(V))$ and exhibit some admissible relations on the object $\End(V)^{\ot 2}$. Below, $L=\|l_i^j\|$ is a matrix and each of the entries $l_i^j$ is identified with the element $x_i\ot x^j\in
V\ot V^*$.

\begin{proposition} The relations{\rm
\be
R^{\varepsilon_1}\, L_1\, R^{\varepsilon_2}\, L_1 -
L_1\, R^{\varepsilon_3}\, L_1\, R^{\varepsilon_4} = 0,
\label{genRE}
\ee}
$\!\!$where $R$ is a Hecke symmetry and $ \varepsilon_i\in \{ 1, -1\}, \, i=1,2,3,4$ are admissible.
\end{proposition}

{\bf Proof.}
In virtue of the identification  $l_i^j = x_i\otimes x^j \in V\otimes V^*$, the coaction (\ref{coac}) leads to the following
transformation for elements of the matrix $L$: $l_i^j\mapsto t_i^ks_p^j\otimes l_k^p$. This coaction can
be symbolically written in the matrix form
\be
L\mapsto T L S,
\label{coa}
\ee
where we assume that the entries of the matrix $L$ commute with those of $S$ and $T$.

According to Definition \ref{def:6} we have to verify, that the linear span, generated by the left hand side\footnote{With
the use of endomorphism $\Psi $ (see definition \ref{psi}) or its analog $\Psi_{R^{-1}}$ for the inverse matrix $R^{-1}$ (see (\ref{psi-r-1}) below) we can write this linear span as the image of the map $\beta = I -\gamma $, where
$\gamma: \End(V)^{\otimes 2}\to \End(V)^{\otimes 2}$ reads
$$
\gamma(L_1L_2) = \Tr_{(0)}\Big(R_{10}^{-\varepsilon_1}L_0R_{10}^{\varepsilon_3}L_0
R_{10}^{\varepsilon_4}(\Psi_{R^{\varepsilon_2}})_{02}\Big)P_{12}.
$$} of (\ref{genRE}) is an invariant subcomodule in $\End(V)^{\otimes 2}$ with respect to the coaction (\ref{coa}). This can
easily be done with the use of the following direct consequence of the RTT relations
$$
S_1\, R^{\varepsilon}\, T_1= T_2\, R^{\varepsilon}\, S_2, \quad \varepsilon\in \{ 1, -1\}.
 $$
Note that the RTT-type relations on the elements $l_i^j$
$$
R_{12}L_1L_2-L_1L_2R_{12} = 0
$$
are not admissible since the corresponding linear span is not invariant comodule with respect to the coaction
(\ref{coa}). So, the RTT algebra is not  braided  (according to  the definition  above). \hfill \rule{6.5pt}{6.5pt}

\begin{definition}{\rm
The algebra $\LR $ generated by matrix element of the matrix $L=\|l_i^j\|$ subject to a particular case of the relations (\ref{genRE})
\be
R\, L_1\, R\, L_1 - L_1\, R\, L_1\, R=0
\label{RE}
\ee
is called the Reflection Equation (RE) algebra. The matrix $L$ is called the generating matrix of the
algebra $\LR$.}
\end{definition}

\begin{remark}{\rm
The relation (\ref{genRE}) with $\varepsilon_1=\varepsilon_4$ and $\varepsilon_2=\varepsilon_3$ is equivalent
to (\ref{RE}) up to a possible change of $R$ for $R^{-1}$. Note that $R^{-1}$ is a braiding. Besides, $R^{-1}$
is a Hecke symmetry if and only if $R$ is also.}
\end{remark}

It was discovered by Sh. Majid that any RE algebra has a braided bialgebra structure (see \cite{M}). This structure is
defined via a braiding $R_{\End(V)}: \End(V)^{\ot 2} \to \End(V)^{\ot 2}$ extended to the whole RE algebra. Note that in
the standard case this braiding is the image of the $U_q(sl(m))$ universal R-matrix in the space $\End(V)^{\ot 2}$.

The relation $R_{\End(V)}(L\otimes L) -L\otimes L =0$ is admissible since it can be rewritten as a particular case of
(\ref{genRE})
$$
R\,L_1R^{-1}L_1 = L_1R\,L_1R^{-1}.
$$
However, in contrast with the RE algebra, the algebra defined by this system of relations does not in general possess
the good deformation property in the sense of \cite{GPS2}.

In this connection we would like to mention the paper \cite{BZ} where another way was suggested of defining analogs of symmetric and skew-symmetric algebras of irreducible modules over the quantum groups. Note that except for the case
$m=2$ the ``R-symmetric algebra" of the adjoint $U_q(sl(m))$-module as defined in \cite{BZ} does not possesses the good
deformation property.

\section{Braided Weyl algebras and their $q\to 1$ limits}

Let $R=R(q)$ be a skew-invertible Hecke symmetry and $\LR $ the corresponding RE algebra defined by relations (\ref{RE})
on their generators. We perform the following linear change of the generators
\be
L=\h\, I-(q-\q)\,\tilde L,\quad \tilde L=\|\tilde l_i^j\|.
\label{shift}
\ee
Then the system (\ref{RE}) being rewritten via the matrix $\tilde L$ becomes
\be
R\, \tilde L_1 R\, \tilde L_1-\tilde L_1 R\, \tilde L_1 R=\h\, (R\, \tilde L_1-\tilde L_1 R).
\label{mRE}
\ee
We call the algebra defined by generators $\hat{l}_i^j$ subject to (\ref{mRE}) the {\em modified} Reflection Equation algebra.
In fact, it is nothing but another form of the RE algebra.

Now, passing to the limit $q\to 1$ we get just the system (\ref{mRE}) but with an involutive braiding $R(1)$. Assuming
$R(1)$ to be the {\it usual} flip $P$ we get the defining relations of the algebra $U(gl(m)_\h)$. If a Hecke symmetry
$R$ is a deformation of the {\it super-flip} acting on the space $\vv_{m,n} $ where $V_{m,n}=V_0\oplus V_1$, $\dim
V_0=m$, $V_1=n$, we get the enveloping algebra $U(gl(m|n)_\h)$ as $q\to 1$. Thus, though the algebras $\Sym(gl(m|n))$
and $U(gl(m|n)_\h)$ are not isomorphic to each other, their $q$-deformations (\ref{RE}) and (\ref{mRE}) are.

In general, as $q\to 1$ we get an algebra which can be considered as the enveloping algebra of a generalized Lie algebra.
Such algebras have been introduced by one of the authors in \cite{G1} (also see \cite{G2}).

Now, we consider an associative unital algebra $\W(\N)$ generated by entries $n_i^j$ and $d_i^j$
of two $\dim V\times \dim V$ matrices $N=\|n_i^j\|$ and $D=\|d_i^j\|$ satisfying the following relations
\be
\begin{array}{rcl}
R\,N_1R\,N_1-N_1R\,N_1R&=&\h\,(R\,N_1-N_1R),\\
R^{-1}D_1R^{-1}D_1&=&D_1R^{-1}D_1R^{-1},\\
D_1R\,N_1R-R\,N_1R^{-1}D_1&=&R+\h\,D_1R.
\end{array}
\label{set}
\ee

We call the algebra $\W(\N)$ the {\em braided Weyl algebra}. Our terminology is motivated by the fact that in the case
$q\to 1$ and $\h\to 0$ (provided $R$ is the standard Hecke symmetry) the algebra $\W(\N)$ is isomorphic
to the usual Weyl algebra generated by commutative indeterminates $n_i^j$ and the partial derivatives
$\partial/\partial n_i^j$ (see (\ref{cl-diff}--\ref{br})). Note that the algebra $\W(\N)$ was introduced in \cite{GPS6}.
It is just the braided Heisenberg double mentioned in the Introduction.

The algebra $\W(\N)$ possesses two subalgebras. The modified RE subalgebra $\N $, generated by the matrix
$N=\|n_i^j\|$, plays the role of a function algebra. The other subalgebra $\D $, generated by the matrix $D=\|d_i^j\|$,
is an RE algebra corresponding to the Hecke symmetry $R^{-1}$. The matrix elements $d_i^j$ are braided analogs
of the partial derivatives on the algebra $\N $. The third set of relations in the system (\ref{set}) (the permutation
relations) plays the role of the Leibniz rule.

In order to show that the algebra defined by (\ref{set}) is a braided algebra indeed, we identify entries of the matrices
$N$ and $D$ with elements of $\End(V)$ as above. Consequently, the coaction of the RTT algebra is given by the
same formula (\ref{coa}) on the both matrices. We leave to the reader checking that the relations (\ref{set}) are
admissible (see Remark \ref{rem7}).

We define the action of the generators $d_i^j$ on the unit $1_\N $ of the subalgebra $\N $ by setting\footnote{Hereafter,
we denote the action of an operator $A$ on an element $x$ in two ways: either $A\triangleright x$ or $A(x)$ depending
on our convenience.} $d_i^j\triangleright 1_\N=0$. Now, the result of action $d_i^j \triangleright n \equiv d_i^j \triangleright
(n \cdot 1_\N)$ on an arbitrary element $n\in \N $ can be found by permuting $d_i^j$ with $n$ and then evaluating it on the
unit element $1_\N $ according to the above rule.

Equivalently, the latter operation can be presented in terms of the counit map $\varepsilon: \D\to {\Bbb K}$ defined by the
following rule
\be
\varepsilon(1_\D)=1,\quad \varepsilon(d_i^j)=0,\quad \varepsilon(d_1d_2)=\varepsilon(d_1)\varepsilon(d_2),
\qquad d_1,d_2\in \D.
\label{coun}
\ee
Checking that the operation $\triangleright$ is well defined on the algebra $\N $ (i.e. the set of defining relations of the algebra $\N $ is
invariant with respect to this action) follows from the results of \cite{GPS5}.

\begin{proposition}
The action of generators $d_i^j$ on the subalgebra $\N $ completed with the natural action of the unit $1_\D\triangleright
n=n$ can be extended up to the action of the whole subalgebra $\D $ on $\N $:
$\D\otimes \N \stackrel{\triangleright}{\rightarrow} \N $. This action gives a representation of the algebra $\D $ in $\N $.
\end{proposition}

This proposition also follows from \cite{GPS5}.

Let now $R$ be an involutive symmetry. Using the permutation relations among $D$ and $N$ we derive a form of the
Leibniz rule convenient for our subsequent considerations.

We introduce the following chains of $R$-matrices:
$$
\RR_{kp} = R_{p-1}R_{p-2}\dots R_{k+1}R_kR_{k+1}\dots R_{p-2}R_{p-1},\quad 1\le k<p.
$$
The Yang-Baxter equation for $R$ allows us to rewrite the above definition in an equivalent form
$$
\RR_{kp} = R_{k}R_{k+1}\dots R_{p-2}R_{p-1}R_{p-2}\dots R_{k+1}R_{k},\quad 1\le k<p.
$$
Besides, we assume that $\RR_{kp} = \RR_{pk}$ by definition.

For any $i,j,k$ we have the following "index exchange rules", well known for the usual flip operators:
$$
\RR_{ij}\RR_{ik} = \RR_{jk}\RR_{ij} = \RR_{ik}\RR_{kj}.
$$

For any braiding $R:\vv\to\vv $ and any $\dim V\times \dim V$ matrix $N$ we introduce a convenient notation:
\be
N_{\overline k}= R_{k-1}\dots R_1N_1R_1^{-1}\dots R_{k-1}^{-1}. \label{conv}
\ee
Then, the following relations are valid
\begin{eqnarray}
&&\RR_{pk}\,N_{\overline k} = N_{\overline p}\, \RR_{pk},\quad k>p, \nonumber\\
\rule{0pt}{5mm}
&&\RR_{kp}\,N_{\overline k} =N_{\overline p} \,\RR_{kp},\quad k<p, \label{n-k-p}\\
\rule{0pt}{5mm}
&&\RR_{pk}\,N_{\overline s} = N_{\overline s}\, \RR_{pk}, \quad \forall\,s\not\in\{p,k\}.\nonumber
\end{eqnarray}
For an involutive symmetry $R$ the relation (\ref{conv}) can be presented as
$$
N_{\overline k} = R_{k-1}\dots R_1N_1R_1\dots R_{k-1}.
$$

Now, we rewrite the permutation relations of $D$ and $N$ (the third line in (\ref{set})) in the form
\be
D_1N_{\overline 2} = N_{\overline 2}D_1 + \hbar D_1R_1+R_1.
\label{D-N}
\ee
It is convenient to modify the generating matrix $D$ as follows
\be
\tilde D = \h^{-1}\,{\rm Id} +D.
\label{d-hat}
\ee
Then, the permutation rules (\ref{D-N}) take the form
$$
\tilde D_1N_{\overline 2} = N_{\overline 2}\tilde D_1 + \hbar\,\tilde D_1\RR_{12}.
$$
The above relation can be easily generalized to
\be
\tilde D_1 N_{\overline k} = N_{\overline k}\tilde D_1 +\hbar\,\tilde D_1 \RR_{1k},
\quad \forall\, k\ge 2.
\label{hatD-N}
\ee
The action of $\tilde D$ on any element of $\N $ can be obtained from the same scheme as above
with the only modification: $\varepsilon({\tilde D})=\h^{-1}$.

Now, with the use of (\ref{hatD-N}) we get the commutation of $\tilde D$ with some low degree
polynomials in $N$
\begin{eqnarray*}
\tilde D_1N_{\overline 2}N_{\overline 3} &=&
N_{\overline 2}N_{\overline 3}\tilde D_1 +\hbar\left(N_{\overline 2}\tilde D_1\RR_{13}+N_{\overline 3}
\tilde D_1\RR_{12} \right) +\hbar^2\tilde D_1\RR_{12}\RR_{23},\\
\rule{0pt}{5mm}
\tilde D_1N_{\overline 2}N_{\overline 3}N_{\overline 4} &=&
N_{\overline 2}N_{\overline 3}N_{\overline 4}\tilde D_1 +\hbar\left(N_{\overline 2}N_{\overline 3}\tilde D_1\RR_{14}+
N_{\overline 2}N_{\overline 4}\tilde D_1\RR_{13} +N_{\overline 3}N_{\overline 4}\tilde D_1\RR_{12}\right)\\
&& +\hbar^2\left(N_{\overline2}\tilde D_1\RR_{13}\RR_{34}+N_{\overline3}\tilde D_1\RR_{12}\RR_{24}
+N_{\overline 4}\tilde D_1\RR_{12}\RR_{23}\right) +\hbar^3\tilde D_1\RR_{12}\RR_{23}\RR_{34}.
\end{eqnarray*}
These relations immediately give us the action of the matrix $D$ on the same polynomials:
\begin{eqnarray*}
D_1\triangleright N_{\overline 2} &=& \RR_{12},\\
\rule{0pt}{5mm}
D_1\triangleright N_{\overline 2}N_{\overline 3} &=&
N_{\overline 2}\RR_{13}+N_{\overline 3}\RR_{12} +\hbar\,\RR_{12}\RR_{23},\\
\rule{0pt}{5mm}
D_1\triangleright N_{\overline 2}N_{\overline 3}N_{\overline 4} &=&
N_{\overline 2}N_{\overline 3}\RR_{14}+ N_{\overline 2}N_{\overline 4}\RR_{13} +
N_{\overline 3}N_{\overline 4} \RR_{12}\\
&& +\hbar\left(N_{\overline2} \RR_{13}\RR_{34}+N_{\overline3} \RR_{12}\RR_{24}
+N_{\overline4} \RR_{12}\RR_{23}\right) +\hbar^2\RR_{12}\RR_{23}\RR_{34}.\\
\end{eqnarray*}

To write down the general result we introduce some more notation. First,
for any set of integers $2\le k_1<k_2<\dots <k_s\le p$ we denote
$$
\RR_{(1k_1\dots k_s)} = \RR_{1k_1}\RR_{k_1k_2}\dots \RR_{k_{s-1}k_s},\quad
$$
an analog of the cycle permutation $(1k_1k_2\dots k_s)$ in the permutation group
$S_p$. Second, for the same set of integers we denote
$$
\left(N_{\overline 2}\dots N_{\overline p}\right)^{(k_1\dots k_s)} =
N_2\dots N_{\overline{k_1-1}}N_{\overline{k_1+1}}\dots
N_{\overline{k_s-1}}N_{\overline{k_s+1}}\dots N_{\overline p},
$$
that is the product $N_{\overline 2}\dots N_{\overline p}$ where the multipliers
$N_{\overline {k_1}}$, $N_{\overline {k_2}}, ..., N_{\overline {k_s}}$ are omitted.

Then, the following result can be proved by induction in $p\ge 2$
\be
\tilde D_1N_{\overline 2}\dots N_{\overline p} =
N_{\overline 2}\dots N_{\overline p}\tilde D_1 +\sum_{s=1}^{p-1}\hbar^s
\sum_{2\le k_1<\dots <k_s\le p}\left(N_{\overline 2}\dots N_{\overline p}\right)
^{(k_1\dots k_s)}\tilde D_1\RR_{(1k_1\dots k_s)}.
\label{D-transp}
\ee
Finally, we get the following formula for the action of the matrix $D$:
\be
D_1\triangleright N_{\overline 2}\dots N_{\overline p} =
\sum_{s=1}^{p-1}\hbar^{s-1}
\sum_{2\le k_1<\dots <k_s\le p}\left(N_{\overline 2}\dots N_{\overline p}\right)
^{(k_1\dots k_s)}\RR_{(1k_1\dots k_s)}. \label{diff}
\ee

Now, we consider the particular case $R=P$, where $P$ is the usual flip. Then, for the first order monomials
formula (\ref{diff}) gives
\be
D_1\triangleright N_2 = P_{12}\qquad {\rm or}\qquad
d_i^j\triangleright n_k^l=P_{ik}^{jl} = \de_i^l\,\de_k^j.
\label{cl-diff}
\ee

First, note that if $\h=0$ (and, therefore, $\N\cong \Sym(gl(m))$) formula (\ref{diff}) assumes that the
action (\ref{cl-diff}) is extended to an arbitrary order monomial via the classical Leibniz rule. This allows us to
identify $d_i^j=\pa/\pa{n_j^i}$.

Second, for $\h\not=0$ we consider the $GL(m)$-covariant product
\be
n_i^j\circ n_k^l=\h\,\de_k^j\, n_i^l,
\label{prod}
\ee
well defined in the algebra $gl(m)_\h $. The Lie bracket in $gl(m)_\h $ is related to this product in the usual
way
\be
[n_i^j, n_k^l]=n_i^j\circ n_k^l-n_k^l\circ n_i^j.
\label{br}
\ee

Now, in the case $R=P$ formula (\ref{diff}) can be obtained by a sequence of several steps. At the first step we apply
an element $d_i^j$ to a given monomial $n_{i_1}^{j_1}\, n_{i_2}^{j_2}... n_{i_p}^{j_p}$ in accordance with the rule
(\ref{cl-diff}) and the classical Leibniz rule. This gives us the sum of $p-1$ terms corresponding to the value $s=1$
of the summation index in (\ref{diff}).

At the second step we choose an arbitrary pair of elements $n_{i_s}^{j_s}$ and $n_{i_t}^{j_t}$, $1\leq s<t\leq p$,
in the initial monomial and calculate their product (\ref{prod}) preserving the order of the elements. After that we
apply $d_i^j$ to the result of this product:
$$
d_i^j\triangleright (n_{i_s}^{j_s}\circ n_{i_t}^{j_t}) = \h\,\delta_i^{j_t}\delta_{i_t}^{j_s}\delta_{i_s}^{j}\qquad{\rm or}
\qquad D_1\triangleright (N_s\circ N_t) = \h\, P_{1s}P_{st}.
$$
This gives us $(p-1)(p-2)/2$ terms proportional to the first order of the parameter $\h$.

In order to get the $\h^2$-order terms we apply the same procedure to all triples of the factors in the given monomial.
Thus, any triple of chosen elements $n_{i_s}^{j_s}$, $n_{i_t}^{j_t}$, $n_{i_r}^{j_r}$ is replaced by the numerical factor
$\h^2\,\de_{i}^{j_r} \, \de_{i_r}^{j_t}\, \de_{i_t}^{j_s}\,\de_{i_s}^{j}$ or $\h^2\,P_{1s}P_{st}P_{tr}$. And
so on: at the $k$-th step we obtain the terms proportional to $\h^{k-1}$.

Note, that in \cite{GPS6} the above procedure of calculating the action of an element $d_i^j$ on a monomial in
$n$ (the ``$\h$-Leibniz rule'') was defined in a different but equivalent way via a coproduct
$\De(d_i^j)=d_i^j\ot 1 +1\ot d_i^j+\h\, d_k^j\ot d_i^k$ which was found by S.~Meljanac and Z.~\v{S}koda.

\begin{remark}
{\rm Note, that the method of computing the action of an element $d_i^j$ on monomials can be readily generalized to
super-algebras $U(gl(m|n)_\h)$. The only difference is that all transpositions of factors in monomials (we have to put
them aside in order to apply the product $\circ $) and their transpositions with elements $d_i^j$ must be done with
taking into account the parity of elements.

A similar version of the Leibniz rule is valid for any "generalized Lie algebra" corresponding to a skew-invertible {\rm
involutive} symmetry $R$ (see \cite{GPS1}). In this case all aforementioned transpositions must be performed in terms
of the involutive symmetry coming in the definition of such a generalized Lie algebra. Unfortunately, if the initial symmetry
$R$ is not involutive (i.e. $q\not=1$), such a simple form of the Leibniz rule is not found.

Our scheme of defining partial derivatives can be also generalized to the enveloping algebras of current Lie algebra
$\widehat{gl(m)}$ (and their super-analogs) since the Lie brackets in these Lie algebras can be realized via a product
$\circ $ analogous to (\ref{prod}). We do not know any way of defining similar partial derivatives on the enveloping
algebras of Lie algebras without such a product, for instance on $U(sl(m))$.}
\end{remark}

\section{Invariant differential operators, an example}

Let $R$ be a skew-invertible Hecke symmetry. It is well known that the elements $\Tr_{\!R}L^k=\Tr(CL^k)$ (the
matrix $C$ is defined in (\ref{BC})), are central in the algebra $\LR $ for any $k$. Note, that a similar statement
is valid for a modified RE algebra, in particular, for the algebra $\N $ coming in the definition (\ref{set}) of the Weyl
algebra $\W(\N)$. Namely, the elements $\Tr_{\!R}N^k$ are central in this algebra. We call the map $L^k\mapsto
\Tr_{\!R}L^k$ {\em  braided trace}.

The subalgebra $\D \subset \W(\N)$ is also an RE algebra with $R$ replaced by $R^{-1}$. For any Hecke symmetry
$R$ with the bi-rank $(m|n)$ we have $R^{-1}=R-(q-q^{-1})I$ and (see \cite{GPS1})
$$
\Tr\,B=\Tr\, C = \frac{(m-n)_q}{q^{m-n}}.
$$

Emphasize that the braiding $R^{-1}$ is also skew-invertible since
the corresponding endomorphism $\Psi_{R^{-1}}$ exists and can be expressed via $\Psi $ (see (\ref{psi})) by the relation
\be
\Psi_{R^{-1}} = \Psi_{12} +(q-q^{-1})q^{2(m-n)}C_1B_2.
\label{psi-r-1}
\ee

This formula enables us to compute the endomorphisms $B_{R^{-1}}$ and $C_{R^{-1}}$. Namely, we have
$$
B_{R^{-1}} = \Tr_{(1)}\Psi_{R^{-1}} = q^{2(m-n)} B,\qquad C_{R^{-1}} = \Tr_{(2)}\Psi_{R^{-1}} = q^{2(m-n)} C.
$$
Thus, the matrices $B_{R^{-1}}$ and $C_{R^{-1}}$ are proportional to these $B$ and $C$ respectively and the elements $\Tr_{\!R} D^k $ differ from central
elements $\Tr_{R^{-1}} D^k= \Tr(C_{R^{-1}}D^k)$ of the subalgebra $\D $ by the factor $q^{-2(m-n)}\not=0$.
So, we associate "braided differential operators" with these elements by treating matrix elements $d_i^j$ of the matrix
$D$ to be derivations on the algebra $\N $ as was described in the previous section. We call $\Tr_{\!R} D^k$
{\em Laplace operators}.

In the same way, we define {\em sl-type Laplace operators} associated with the elements $\Tr_R \tilde{D}^k$ where
\be
\tilde{D}=D-{{\Tr_{\!R} D}\over{\Tr_{\!R} I}}\,I
\label{less}
\ee
is the R-traceless part of the matrix $D$. Since
$$
\Tr_{\!R} I=\Tr C=q^{n-m}(m-n)_q,
$$
the operator $\tilde{D}$ is well defined provided $m\not=n$ and $q$ is subject to the condition (\ref{cond}). Note,
that the families of the operators $D^k$ and $\tilde{D}^k$ can be expressed via each other.

Also, we need the operators $\Tr_{\!R}( N^lD^k)$, $\Tr_{\!R}(N^l\tilde{D}^k)$  and $\Tr_{\!R}(\tilde N^l\tilde{D}^k)$
which are called {\em invariant}. In the
standard case they are invariant with respect to the action of the quantum group $\Uq$. In general, they are invariant with
respect to the coaction of the corresponding RTT algebra.

We are interested in the action of the Laplace operators on the algebra $\N$, in particular, on its center
$Z(\N)$.

\begin{conjecture} All the Laplace operators $\Tr_{\!R}D^k$ map the center $Z(\N)$ into itself.
\end{conjecture}

Below, we consider invariant differential operators in a particular case $R=P$, $m=2$, $n=0$. Also we assume
${\Bbb K} = {\Bbb R}$. The case ${\Bbb K} = {\Bbb C}$ can be treated as the complexification of the objects and morphisms over the real field.

Let us recall some results and notations from \cite{GPS6}. Consider the algebra $U(u(2)_\h)$ with the following
defining relations between generators
$$
[x, \, y]=\h\, z,\quad [y, \, z]=\h\, x,\quad[z, \, x]=\h\, y,\quad [t, \, x]=[t, \, y]=[t, \, z]=0.
$$
These generators are connected with the generators $a$, $b$, $c$ and $d$ of the algebra $U(gl(2)_\h)$
in the standard way
\be
a=t-iz,\qquad b=-ix-y,\qquad c=-ix+y,\qquad d = t+iz.
\label{abcd}
\ee

The center $Z(U(u(2)_\h))$ is generated by $t$ and the Casimir element
$$
\Cas = x^2+y^2+z^2.
$$

First, we consider the algebra $U(su(2)_\h)$ which is known to be flat over its center (see \cite{D}).
More precisely, as an $su(2)$-module it is isomorphic to the following direct sum of $su(2)$-modules
\be
U(su(2)_\h)\cong \bigoplus_{k =0}^\infty\Bigl(Z(U(su(2)_\h))\otimes V^k\Bigr),
\label{sum}
\ee
where $V^k$ is the $su(2)$-module of the highest weight $k$. We fix its highest weight vector
in the complexification of the algebra $U(su(2)_\h)$ to be $b^k=(-i\, x-y)^k$.
\begin{remark}
\label{rem2}
{\rm
In fact, the module $V^k$ contains the elements Re$(b^k)$ and Im$(b^k)$.
}
\end{remark}

\begin{definition} {\em The component $Z(U(su(2)_\h))\otimes V^k$ is called isotypic.}

\end{definition}

The algebra $U(u(2)_\h)$ is the first ingredient of the Weyl algebra we are going to construct. The second
ingredient is the algebra $\D$ generated by the patrial derivatives $\pa_x$, $\pa_y$, $\pa_z$ and $\tilde \pa_t
= \pa_t + \frac{2}{\h}\,I$. This shift of the $t$-derivative is convenient for technical reasons.

The properties of the partial derivatives follow from the defining relations (\ref{set}). Namely, one can show
(see \cite{GPS6} for detail) that the partial derivatives commute with each other, while their permutation relations
 with the generators $x$, $y$, $z$ and $t$ read
\be
\begin{array}{l@{\quad}l@{\quad}l@{\quad}l}
\tilde\pa_t\,t - t\,\tilde\pa_t = \hh\,\tilde\pa_t & \tilde\pa_t\, x - x\,\tilde\pa_t
=-\hh\,\pa_x &
\tilde\pa_t\, y - y\, \tilde\pa_t=-\hh\,\pa_y &\tilde\pa_t\, z - z\,\tilde\pa_t=- \hh\,\pa_z\\
\rule{0pt}{7mm}
\pa_x\, t - t\,\pa_x = \hh\,\pa_x &\pa_x \,x - x\,\pa_x = \hh\,\tilde\pa_t &
\pa_x \, y- y\,\pa_x = \hh\,\pa_z & \pa_x \,z - z\, \pa_x = - \hh\,\pa_y \\
\rule{0pt}{7mm}
\pa_y \,t - t \, \pa_y = \hh\,\pa_y & \pa_y \,x - x\, \pa_y = -\hh\,\pa_z &
\pa_y \,y - y \, \pa_y = \hh\,\tilde\pa_t & \pa_y \,z - z \, \pa_ y= \hh\,\pa_x\\
\rule{0pt}{7mm}
\pa_z \,t - t \,\pa_z = \hh\,\pa_z & \pa_z \,x - x \,\pa_z = \hh\,\pa_y&
\pa_z \,y - y\,\pa_z = -\hh\,\pa_x & \pa_z \,z - z \,\pa_z = \hh\,\tilde\pa_t.
\end{array}
\label{leib-r}
\ee

Our next aim is to compute the action of the invariant operators $\tpa$ and $\Tr_{\!R}\tilde D^2$ on elements
$t^k $ and $f b^k$, $\forall f\in Z(U(su(2)_\h))$ from the complexification of the algebra $U(u(2)_\h)$. Since
the action of these operators is $SU(2)$-invariant (or $SL(2)$-invariant if we consider the complexification), we
easily extend this action up to the isotypic component $Z(U(su(2)_\h)) \otimes V^k$ and, consequently, on the whole
algebra $U(u(2)_\h)$ due to the isomorphism (\ref{sum}). Note, that on the algebra $U(su(2)_\h)$ the Laplace
operator $\Tr_{\!R}\tilde D^2$ coincides (up to a numerical factor) with the invariant differential operator
$$
\Delta = \pa_x^2+\pa_y^2+\pa_z^2.
$$
Below we deal with $\Delta$ instead of $\Tr_{\!R}\tilde D^2$.

In contrast with the classical case $\h=0$, when the function algebra is $\Sym(u(2))$, the computation of the
action of an invariant operator  on $U(u(2)_\h)$ involves other invariant operators. Bellow, we will deal with two invariant first
order differential operators
\be
\tilde\partial_t \quad {\rm and} \quad Q =x\, \pa_x+y\,\pa_y+z\,\pa_z
\label{first}
\ee
and with four invariant second order differential operators
\be
\De_0=\dd^2,\qquad \De_1=\De=\pa_x^2+\pa_y^2+\pa_z^2,\qquad \De_2=Q\,\dd,\qquad \De_3=Q^2.
\label{second}
\ee
It is evident, that any polynomial in the above operators with coefficients from $Z(U(u(2)_\hbar))$ is also an
invariant operator.

First, we find the permutation relations of the operators (\ref{first}) and (\ref{second}) with the central elements $t^k$,
$k\in {\Bbb N}$:
\be
\tpa \, t^k= (t+{\h\over 2})^k\, \tpa, \qquad
Q \, t^k= (t+{\h\over 2})^k\,Q,\qquad
\De_i\, t^k= (t+\h)^k\, \De_i,\quad \forall\, i=0,1,2,3.
\label{act-t}
\ee

Second, the permutation relations of the first order operators (\ref{first}) with the Casimir element $\Cas$
are as follows
\begin{eqnarray*}
&&\tilde\partial_t\,\Cas =\left(\Cas-\frac{3}{4}\,\h^2\right)\tilde\partial_t-\h\,Q,\\
\rule{0pt}{7mm}
&& Q\,\Cas =\h\,\Cas\,\tilde\partial_t+ \left(\Cas+\frac{\h^2}{4}\right)Q,
\end{eqnarray*}
or in the matrix form
\be
\left(
\begin{array}{c}
 \tilde\partial_t\\
\rule{0pt}{6mm}
Q
\end{array}
\right)\Cas =
\left(\begin{array}{c@{\hspace*{6.5mm}}c}
\displaystyle \Cas -{3 \h^2 \over 4}& -\h\ \\
\rule{0pt}{3mm}
\h\Cas &\displaystyle \Cas+{\h^2 \over 4}
\end{array}\right)
\left(
\begin{array}{c}
 \tilde\partial_t\\
\rule{0pt}{6mm}
Q
\end{array}
\right).
\label{cas-phi}
\ee
We denote $\Phi=\Phi(\Cas)$ the $2\times 2$ matrix coming in this formula. Then the permutation relations of the
operators (\ref{first}) with the central element $\Cas^p,\, p\in \NN$ can be expressed via the $p$-th power of the matrix
$\Phi$.

Now, recall that the generating matrix
\be N=\left(\!\begin{array}{cc}
t-i\,z& -i\, x-y\\
\rule{0pt}{5mm}
-i\, x+y& t+i\, z
\end{array}\!\right)
\label{matrN}
\ee
of the algebra $U(u(2)_\h)$ satisfies the Cayley-Hamilton (CH) identity of the form
$$
N^2-(2\,t+\h)\, N+\,(t^2+x^2+y^2+z^2+\h\, t)\,I= 0.
$$
Let $\mu_1$ and $\mu_2$ obey the relations
$$
\mu_1+\mu_2=2\,t+\h,\qquad
\mu_1 \,\mu_2= t^2+x^2+y^2+z^2+\h\, t.
$$
This means that the variables $\mu_1$ and $\mu_2$ belong to an algebraic extension of the center $Z(U(u(2)_\h))$.
We call them the {\em eigenvalues} of the matrix $N$. Note, that our ordering of the eigenvalues is arbitrary.

These eigenvalues are useful for parameterizing all central elements. In particular, for the Casimir element $\Cas$
we get
\be
\Cas=\frac{\h^2-\mu^2}{4},\quad \mu = \mu_1-\mu_2.
\label{cas-mu}
\ee

Also, we need the eigenvalues $\lambda_1$ and $\lambda_2$ of the matrix $\Phi(\Cas)$
$$
\lambda_1 = \frac{ \mu}{4}\,(2\h - \mu),\qquad
\lambda_2 = - \frac{ \mu}{4}\,(2\h + \mu).
$$
The matrix $\Phi(\Cas)$ has the following spectral decomposition
$$
\Phi(\Cas) = \lambda_1\cdot \frac{\Phi-\lambda_2}{\lambda_1-\lambda_2}
+\lambda_2\cdot \frac{\Phi-\lambda_1}{\lambda_2-\lambda_1} = \lambda_1P_1(\Phi)
+\lambda_2 P_2(\Phi),
$$
where the explicit form of the matrices $P_i(\Phi)$ is as follows
$$
P_1(\Phi) =
\frac{1}{ \mu}
\left(\begin{array}{c@{\hspace*{6.5mm}}c}
\displaystyle { \mu-\h \over 2}&\displaystyle -1 \\
\rule{0pt}{8mm}
\displaystyle {\h^2-\mu^2\over 4} &\displaystyle { \mu+\h \over 2}
\end{array}\right),\qquad
P_2(\Phi) = \frac{1}{ \mu}
\left(\begin{array}{c@{\hspace*{6.5mm}}c}
\displaystyle { \mu+\h \over 2}&\displaystyle 1 \\
\rule{0pt}{8mm}
\displaystyle {\mu^2 -\h^2\over 4} &\displaystyle { \mu-\h \over 2}
\end{array}\right).
$$
The matrices $P_i(\Phi)$, $i=1,2$, are complementary projectors, i.e. they satisfy the relations
$$
P_i(\Phi)P_j(\Phi) = \de_{ij}P_i(\Phi),\qquad P_1(\Phi)+P_2(\Phi) =I.
$$
The spectral decomposition enables us to compute the matrix $\Phi^p(\Cas)$ explicitly:
\begin{eqnarray}
\Phi^p= \lambda_1^pP_1(\Phi)+\lambda_2^pP_2(\Phi)=&& \nonumber\\
&&\hspace*{-20mm} \rule{0pt}{12mm}\frac{\la_1^p}{ \mu}
\left(\begin{array}{c@{\hspace*{6.5mm}}c}
\displaystyle { \mu-\h \over 2}&\displaystyle -1 \\
\rule{0pt}{8mm}
\displaystyle {\h^2-\mu^2\over 4} &\displaystyle { \mu+\h \over 2}
\end{array}\right)+
\frac{\la_2^p}{ \mu}
\left(\begin{array}{c@{\hspace*{6.5mm}}c}
\displaystyle { \mu+\h \over 2}&\displaystyle 1 \\
\rule{0pt}{8mm}
\displaystyle {\mu^2 -\h^2\over 4} &\displaystyle { \mu-\h \over 2}
\end{array}\right).
\label{one}
\end{eqnarray}

Now, we are able to compute the action of the operators (\ref{first}) on the elements $\Cas^p\, b^k\,\, p,k\in\NN$.
To this end we first apply these operators to the elements $b^k$.

\begin{proposition}
\label{prop:14}
The following relations hold true
\be
\tpa(b^k)={{2}\over{\h}} b^k,\qquad Q(b^k)=kb^k
\label{act}
\ee
and, consequently, $\pa_t(b^k)=0$.
\end{proposition}
\medskip

\noindent{\bf Proof.} To prove formulae (\ref{act}) it is sufficient to use the $\h$-Leibniz rule, described at the end of the
previous section. Namely, since $b\circ b = 0$ (see (\ref{prod})), the action of derivatives on $b$ is, actually, classical
(the terms proportional to $\h$ are absent). This immediately involves the result (\ref{act}).\hfill\rule{6.5pt}{6.5pt}
\medskip

Thus, in order to compute the action of $\tilde\pa_t$ and $Q$ on $\Cas^p\, b^k$ we transpose these operators
with the element $\Cas^p$ with the help of (\ref{cas-phi}) and then apply them to the element $b^k$ in accordance
with (\ref{act}). Finally, we get
$$
\left(
\begin{array}{c}
 \tilde\partial_t\\
\rule{0pt}{6mm}
Q
\end{array}
\right)(\Cas^p \, b^k) = \Phi^p \left(
\begin{array}{c}\displaystyle
{2\over \h}\\
\rule{0pt}{6mm}
k
\end{array}\right)\,b^k=b^k \left(
\begin{array}{c}
\displaystyle {1\over \h}(\la_1^p+\la_2^p)-{(\la_1^p-\la_2^p) \over \mu}(k+1)\\
\rule{0pt}{8mm}\displaystyle
(\la_1^p+\la_2^p){k\over 2}+{(\la_1^p-\la_2^p) \over \mu}\left({\h (k+1)\over 2}-{\mu^2 \over 2\h}\right)
\end{array}\right)
$$
or, in more detail,
\be
\begin{array}{l}
\displaystyle
\pa_t(\Cas^p \, b^k)=b^k\left({1\over \h}(\la_1^p+\la_2^p-2)-{(\la_1^p-\la_2^p) \over \mu}(k+1) \right),\\
\rule{0pt}{8mm}\displaystyle
Q(\Cas^p \, b^k) = b^k\left((\la_1^p+\la_2^p){k\over 2}+{(\la_1^p-\la_2^p) \over \mu}\left({\h (k+1)\over 2}-
{\mu^2 \over 2\h}\right) \right).
\end{array}
\label{acti}
\ee

These formulae enable us to compute the action of the operators (\ref{first}) on an element $f(\Cas)\, b^k$
where $f(u)$ is a polynomial or a formal series in one variable. This action is defined by formula (\ref{acti})
where $\la_1^p$ and $\la_2^p$ are respectively replaced by $f(\la_1)$ and $f(\la_2)$.

However, we are interested in the action of the operators (\ref{first}) on $\mu$ and its integer powers.
To find this action, we take into account that $\mu=\sqrt{\rule{0pt}{3.4mm}\h^2-4\Cas}$ and choose
$f(u)=\left(\!\sqrt{\rule{0pt}{3.4mm}\h^2-4 u}\right)^p$, where $p\in \Z$ is an integer (may be negative). Then we get
$$
f(\la_1)=\left(\sqrt{(\mu-\h)^2}\right)^p=(\mu-\h)^p,\qquad f(\la_2)=\left(\sqrt{(\mu+\h)^2}\right)^p
=(\mu+\h)^p.
$$
Our choice of the roots in these formulae is motivated by the classical limit $\h\to 0$. Finally, we get the
following proposition.

\begin{proposition}
\label{th:15}
The following relations take place (here $k\in \NN$ and $p\in \Z$)
$$
\pa_t(\mu^p \, b^k)=b^k\left({1\over \h}\Bigl((\mu-\h)^p+(\mu+\h)^p-2\Bigr)-
{(\mu-\h)^p-(\mu+\h)^p \over \mu}(k+1) \right),
$$
$$
Q(\mu^p \, b^k)=b^k\left(\Bigl((\mu-\h)^p+(\mu+\h)^p\Bigr){k\over 2}+{(\mu-\h)^p-
(\mu+\h)^p \over \mu}\left({\h (k+1)\over 2}-{\mu^2 \over 2\h}\right)\right).
$$
\end{proposition}

For any rational function $f(u)$ the action of the operators $\pa_t$ and $Q$ on elements $f(\mu)\, b^k$ can be
generalized in the obvious way: we replace $(\mu-\h)^p$ and $(\mu+\h)^p$ in the right hand side of the above formulae
by $f(\mu-\h)$ and $f(\mu+\h)$ respectively. Consequently, by assuming $f(t, \mu)$ to be a rational function in two
variables we define the action of the operators $\pa_t$ and $Q$ on elements $f(t, \mu) b^k$ by similar formulae but with
$t$ in the right hand side replaced by $t+\frac{\h}{2}$.

\begin{remark}{\rm
As follows from the proposition \ref{th:15}, the operators $b^{-k} \pa_t b^k$ and $b^{-k} Q b^k$ (gauge equivalent to $\pa_t$
and $Q$ respectively) are well defined on the space of polynomials (series) in $\mu$. Also, observe that the factor $b^k$ in
all formulae above can be replaced by any element from the module $V^k$.}
\end{remark}

Now, we turn to computing the action of the second order operators (\ref{second}) on the elements
$\Cas^p\, b^k$. To this end we reproduce the permutation relations of these operators with the Casimir
element $\Cas$.

\begin{proposition} {\bf \cite{GPS6}}
The following permutation relations hold true:
\be
\Delta_i \Cas = \sum_{j=0}^3\Pi_{ij}\,\Delta_j, \quad 0\le i\le 3, \label{two}
\ee
where the matrix $\Pi = \Pi(\Cas)$ reads
$$
\Pi =
\left(\begin{array}{c@{\hspace*{6.5mm}}c@{\hspace*{6.5mm}}
c@{\hspace*{6.5mm}}c}
\displaystyle\Cas-{3\over 2}\,\h^2 &\displaystyle{\h^2\over 2} &-2\h &0 \\
\rule{0pt}{8mm}
\displaystyle {3\over 2}\,\h^2 &\displaystyle \Cas-{\h^2\over 2} &2\h& 0\\
\rule{0pt}{8mm}
\h\,\Cas &0 &\displaystyle\Cas-{\h^2\over 2} &-\h \\
\rule{0pt}{8mm}
\h^2\Cas&\displaystyle -{\h^2\over 2}\,\Cas &
\displaystyle\h\left(2\Cas+{\h^2\over 4}\right)&
\displaystyle\Cas+{\h^2\over 2}
\end{array}\right).
$$
\end{proposition}

The matrix $\Pi $ plays the same role for the operators (\ref{second}) as $\Phi $ plays for the operators (\ref{first}).
Thus, we get
$$
\Delta_i\Cas^p = \sum_{j=0}^3(\Pi^p)_{ij}\Delta_j,\,\, p\in \NN.
$$

By a direct calculation one can verify that the matrix $\Pi $ is semisimple:
$$
\Pi\sim {\rm diag}(\lambda_0,\lambda_0,\lambda_+,\lambda_-),
$$
where (recall that $\mu = \mu_1-\mu_2$)
$$
\lambda_0= \frac{1}{4}\,(\h^2 - \mu^2),\quad
\lambda_\pm = \frac{1}{4}\, (\h^2-(\mu \pm 2\,\h)^2)\,.
$$

Similarly to (\ref{one}) we get
\be
\Pi^p = \sum_{a=0,\pm}\lambda_a^p\prod_{b=0,\pm\atop b\not=a}
\frac{(\Pi-\lambda_bI)}{(\lambda_a-\lambda_b)}=
\lambda_0^p P_0(\Pi)+\lambda_+^p P_+(\Pi)+\lambda_-^p P_-(\Pi).
\label{LS-pol}
\ee
Here $P_a(\Pi)$, $a\in\{0,+,-\}$, are complementary projectors:
$$
P_a(\Pi) P_b(\Pi)=\delta_{ab}P_a(\Pi),\qquad
P_0(\Pi)+P_+(\Pi) +P_-(\Pi)=I.
$$
Their explicit form reads:
$$
P_0(\Pi) =
\left(\begin{array}{c@{\hspace*{6.5mm}}c@{\hspace*{6.5mm}}
c@{\hspace*{6.5mm}}c}
\displaystyle {1 \over 2} & 0 & \displaystyle{\h\over {\h^2-\mu^2}} &
 \displaystyle{2\over {\h^2-\mu^2}}\\
\rule{0pt}{8mm}
\displaystyle {1 \over 2} & 1 &\displaystyle{-\h\over {\h^2-\mu^2}} &
 \displaystyle{-2\over {\h^2-\mu^2}}\\
\rule{0pt}{8mm}
\displaystyle {-\h \over 4} & \displaystyle {\h \over 4} & \displaystyle{-\h^2\over {\h^2-\mu^2}} &
 \displaystyle{-2\h\over {\h^2-\mu^2}}\\
\rule{0pt}{8mm}
\displaystyle{{2\h^2-\mu^2}\over 8} &\displaystyle{-\h^2\over 8}&
\displaystyle{{\h(3\h^2-\mu^2)}\over{4(\h^2-\mu^2)}}&
\displaystyle{{3\h^2-\mu^2}\over{2(\h^2-\mu^2)}}
\end{array}\right)
$$

$$
P_+(\Pi) =
\left(\begin{array}{c@{\hspace*{6.5mm}}c@{\hspace*{6.5mm}}
c@{\hspace*{6.5mm}}c}
\displaystyle {2\h+\mu \over 4\mu} &\displaystyle {-\h \over 4\mu}
& \displaystyle{{3\over 2}\h +\mu \over \mu(\h +\mu)} &
 \displaystyle{1 \over \mu(\h +\mu)}\\
\rule{0pt}{8mm}
-\,\displaystyle {2\h+\mu \over 4\mu} &\displaystyle {\h \over 4\mu}
 &-\,\displaystyle{{3\over 2}\h +\mu \over \mu(\h +\mu)} &
 \displaystyle{-1 \over \mu(\h +\mu)}\\
\rule{0pt}{8mm}
\displaystyle{-2\h^2+\h\,\mu+\mu^2\over 8 \mu} &
\displaystyle{\h\,(\h-\mu)\over 8\mu} &
\displaystyle{-3\h^2+\h\,\mu+2\mu^2\over 4 \mu( \h+\mu)}
&\displaystyle{-\h+\mu \over 2\mu(\h+\mu)}\\
\rule{0pt}{8mm}
\displaystyle{-2\h^2+\h\,\mu+\mu^2\over 16} &\displaystyle{\h\,(\h-\mu)\over 16}
& \displaystyle{-3\h^2+\h\,\mu+2\mu^2\over 8( \h+\mu)}&
\displaystyle{-\h+\mu \over 4(\h+\mu)}
\end{array}\right)
$$

$$
P_-(\Pi) =
\left(\begin{array}{c@{\hspace*{6.5mm}}c@{\hspace*{6.5mm}}
c@{\hspace*{6.5mm}}c}
\displaystyle {- 2\h+\mu \over 4\mu} &\displaystyle {\h \over 4\mu}
&\displaystyle{{3\over 2}\h -\mu \over \mu(\mu - \h)} &
 \displaystyle{1 \over \mu(\mu-\h)}\\
\rule{0pt}{8mm}
\displaystyle {2\h - \mu \over 4\mu} &\displaystyle { -\h \over 4\mu}
&\displaystyle{- {3\over 2}\h + \mu \over \mu(\mu - \h)} &
 \displaystyle{ -1 \over \mu(\mu-\h)}\\
\rule{0pt}{8mm}
\displaystyle{2\h^2+\h\,\mu - \mu^2\over 8 \mu} &
\displaystyle{-\h\,(\h+\mu)\over 8\mu} &
\displaystyle{3\h^2+\h\,\mu-2\mu^2\over 4 \mu( \h -\mu)}
&\displaystyle{\h+\mu \over 2\mu(\h-\mu)}\\
\rule{0pt}{8mm}
\displaystyle{-2\h^2-\h\,\mu+\mu^2\over 16} &\displaystyle{\h\,(\h+\mu)\over 16}
& \displaystyle{-3\h^2-\h\,\mu+2\mu^2\over 8( \h-\mu)}&
-\,\displaystyle{\h+\mu \over 4(\h-\mu)}
\end{array}\right)
$$

Now, compute the action of the operators $\De_i$ on the elements $b^k$. Applying the same reasoning as in the
proposition \ref{prop:14}, we get that the $\h $-Leibniz rule reduces to the classical one. Thus, we come to the proposition below.
\begin{proposition}
The following action holds true
$$
\De_0(b^k)={{4}\over{\h}^2} b^k,\quad \De_1(b^k)=0,\quad \De_2(b^k)={{2k}\over{\h}}\, b^k,\quad
\De_3(b^k)=k^2\, b^k.
$$
\end{proposition}

Thus, the action of the operators $\pa_t,\,Q,\,\De $ on the elements $b^k$ and, consequently, on the whole modules $V^k$
are given by the same formulae as in the classical case. Therefore, the elements from $V^k$ are {\em harmonic} in the usual
sense of the words: whey are killed by  the operators $\tpa$ and $\De$.

Using the same method as above, we can compute the action of the operators (\ref{second}) on elements $\Cas^p\, b^k$. Namely, we have
$$
\left(\begin{array}{c}
\De_0\\
\De_1\\
\De_2\\
\De_3
\end{array}\right)(\Cas^p\, b^k)=b^k\Bigl(\la_0^p\, P_0(\Pi)+\la_+^p\, P_+(\Pi)+\la_-^p \,P_-(\Pi)\Bigr)
\left(\begin{array}{c}
\displaystyle{{4}\over{\h^2}}\\
\rule{0pt}{6mm}
0\\
\rule{0pt}{6mm}
\displaystyle
{{2k}\over{\h}}\\
\rule{0pt}{6mm}
k^2
\end{array}\right).
$$

Now, we want to extend this formula to elements $\mu^p\, b^k,\,\,p\in \Z, k\in \NN$. Since
$$
\sqrt{\rule{0pt}{3.4mm}\h^2-4\la_0} = \mu,\qquad \sqrt{\h^2-4\la_\pm} =\mu\pm 2\h
$$
(as above the signs are motivated by the classical limit) we have
\be
\left(\begin{array}{c}
\De_0\\
\De_1\\
\De_2\\
\De_3
\end{array}\right)(\mu^p\, b^k)=b^k\Bigr(\mu^p P_0(\Pi)+(\mu+2\h)^p P_+(\Pi)+(\mu-2\h)^p P_-(\Pi)\Bigl)
\left(\begin{array}{c}
\displaystyle
{{4}\over{\h^2}}\\
\rule{0pt}{6mm}
0\\
\rule{0pt}{6mm}
\displaystyle
{{2k}\over{\h}}\\
\rule{0pt}{6mm}
k^2
\end{array}\right). \label{sec} \ee

This enables us to compute the action of the operators in question on elements $f(\mu)b^k$ where $f(u)$ is a rational function.
 For example, the action of the operator $\De_1=\De$ in the case $k=0$ is as follows:
\be
\De(f(\mu))={{1}\over \h^2}\Bigr(2f(\mu)-f(\mu-2\h)-f(\mu+2\h)\Bigl)+{{2}\over \mu \h}\,\Bigl(f(\mu-2\h)-
f(\mu+2\h)\Bigr).
\label{mmu}
\ee
This formula together with (\ref{act-t}) entails that if $f$ is a rational function in two variables then
\begin{eqnarray}
\De(f(t, \mu))={{1}\over \h^2}\Bigl(2f(t+\h, \mu)-f(t\!\!\!\!&+&\!\!\!\!\h, \mu-2\h)-f(t+\h, \mu+2\h)\Bigr)
\nonumber\\
&+&\!\!\!\!\displaystyle{{2}\over \mu \h}\,\Bigl(f(t+\h, \mu-2\h)-f(t+\h, \mu+2\h)\Bigr).
\label{two-f}
\end{eqnarray}

\begin{remark}{\rm
Formula (\ref{mmu}) (expressed via $(\mu_1-\mu_2)^2$) has been found in \cite{GPS6}. The currently used
element $\mu=\mu_1-\mu_2$ is much more convenient since the expression $\sqrt{\rule{0pt}{3.4mm}\h^2-4u}$,
$u\in\{\la_1,\la_2,\, \la_0,\, \la_\pm\}$ is a first order polynomial in $\mu $. This property considerably simplifies all
calculations.}
\end{remark}

Other formulae expressing the action of the operators (\ref{second}) on elements $f(\mu)b^k$ (and consequently,
on elements $f(\mu)u$ with any $u\in V^k$) can be deduced from relation (\ref{sec}) in the same way. However, they
are too cumbersome and we do not exhibit their explicit form. Finally, the following proposition takes place.

\begin{proposition}
\label{prop:20}
Each of the operators (\ref{first}) and (\ref{second}), maps elements of the form $f(\mu)u$, where $f(\mu)$ is a rational
function and $u\in V^k$, into elements of the same form, namely, $g(\mu)u$ with a rational $g$. The same claim is valid if
$f(t,\mu)$ is a ration function in two variables.
\end{proposition}

\section{Quantization of differential operators and dynamical models}

In this section we use different notation for classical and quantum objects. Namely, for elements of the algebra
$U(u(2)_\h)$, as well as for derivatives on this algebra, we use the Latin letters with hats (for instance, $\hat x$,
$\hat\pa_x$ etc.), while for elements of the commutative algebra $\Sym(u(2))$ and the corresponding derivatives
we keep the previous notation ($x$, $\pa_x$ etc.).

Given a differential operator $\D $ on the algebra $\Sym(u(2))$, we shall associate with it an operator acting on the
algebra $U(u(2)_\h)$. (We are mainly interested in operators describing dynamical models.)

If such an operator has constant (i.e., numerical) coefficients, we define its analog on $U(u(2)_\h)$ by the same formula
but with a new meaning of the derivatives (i.e. all partial derivatives are assumed to act on the algebra $U(u(2)_\h)$
and to be subject to (\ref{leib-r})). Thus, the d'Alembert operator
$$
\hat\square=\hat\pa_t^2-\hat\pa_x^2-\hat\pa_y^2-\hat\pa_z^2
$$
is well defined on the algebra $U(u(2)_\h)$.

In a similar manner we can define quantum analogs of the Dirac and Maxwell operators. We introduce the quantum
analog of the Dirac operator by the classical formula employing the usual Dirac matrices:
$$
\hat D = \gamma^0\hat\partial_t -\gamma^1\hat\partial_x-\gamma^2\hat\partial_y-\gamma^3\hat\partial_z.
$$
As usual,  the Dirac matrices $\gamma^\mu $ realize a representation of the Clifford algebra
$$
\gamma^\mu\gamma^\nu +\gamma^\nu\gamma^\mu = 2g^{\mu\nu}I,\qquad
\mu,\nu\in\{0,1,2,3\}
$$
with the standard Minkowski metric $g^{\mu\nu} = \diag (1,-1,-1,-1)$. Since the {\em quantum partial derivatives}
$\hat\partial_t,...,\hat\partial_x$ commute with each other, the relation between the quantum d'Alembert and
Dirac operators does not change compared with the classical case
$$
\hat D^2 = \hat\square\, I.
$$

For defining the quantum analog of the Maxwell operator we first identify a differential 1-form\footnote{In Electrodynamics
the components of this form are those of the four-potential $A_{\mu}$, $\mu\in\{0,1,2,3\}$. The Maxwell operator comes in
the left hand side of the well-known differential equation $\square A_\mu - \pa_\mu(\pa\cdot A) = 0$.}
$$
\omega=\al\, dt+\beta\, dx+\gamma\, dy+\delta\,dz,\qquad
\al,\,\beta,\,\gamma,\,\delta\in \Sym(u(2))
$$
with the vector-function $(\al,\,\beta,\,\gamma,\,\delta)^t$ (here, the upper script $t$ stands for the transposing)
and realize the classical Maxwell operator \mbox{\sf Mw} as follows
$$
\mbox{\sf Mw}
\left(\!\!
\begin{array}{c}
\al\\
\beta\\
\gamma\\
\delta
\end{array}\!\!
\right)=
\left(\!\!
\begin{array}{c}
\square(\al)\\
\square(\beta)\\
\square(\gamma)\\
\square(\delta)
\end{array}\!\!
\right)-
\left(\!\!
\begin{array}{c}
\partial_t\\
\partial_x\\
\partial_y\\
\partial_z
\end{array}\!\!
\right)\,(\partial_t\,-\partial_x,\,-\partial_y,\,-\partial_z)
\left(\!\!
\begin{array}{c}
\al\\
\beta\\
\gamma\\
\delta
\end{array}\!\!
\right).
$$
Its quantum analog acts on a vector-function $(\hat\al,\,\hat\beta,\,\hat\gamma,\,\hat\delta)^t$, where $\hat\al $,
$\hat \beta $, $\hat\gamma $, $\hat\delta $ are some elements of $U(u(2)_h)$, by the same formula with the change
of the classical partial derivatives for the quantum ones: $\pa\to\hat\pa$, $\square\to\hat \square$.

If the coefficients of a given differential operator $\D $ are some non-numerical elements of $\Sym(u(2))$ the problem of
its quantization becomes more subtle since it involves a quantization of the function algebra $\Sym(u(2))$. Namely, we
need an $SU(2)$-isomorphism of linear spaces $\al: \Sym (u(2))\to U(u(2)_\h)$ satisfying the following properties.

The induced product
$$
f\star_\h g=\al^{-1}(\al(f)* \al(g)),\qquad f,g\in \Sym(u(2))
$$
in the algebra $\Sym(u(2))$ should smoothly depend on the parameter $\h $ (here $*$ stands for the product in
the algebra $U(u(2)_\h)$) so that
$$
\lim_{\h\to 0} f\star_h g = f\cdot g \,\,\,{\rm and}\,\,\, \lim_{\h\to 0}\frac{f\star_\h g - g\star_\h f}{\h} = \{f,g\}_{u(2)}.
$$
Here $f\cdot g$ is the product in the algebra $\Sym(u(2))$ and the notation $\{\,,\,\}_{u(2)}$ stands for the linear
Poisson-Lie bracket on the algebra $\Sym(u(2))$.

Besides, we want the restriction of $\al $ to the Poisson center\footnote{By the Poisson center we mean the set of
elements from $\Sym(u(2))$ central with respect to the Poisson-Lie brackets
$$
\{x,y\}_{u(2)}=z,\quad \{y,z\}_{u(2)}=x,\quad \{z,x\}_{u(2)}=y,\quad \{t,g\}_{u(2)}=0,\quad g\in\{x,y,z\}.
$$
This center is generated by the elements $t$ and $\cas=x^2+y^2+z^2$.} $Z(\Sym (u(2)))$ of the algebra $\Sym (u(2))$,
to be an algebraic isomorphism between $Z(\Sym (u(2)))$ and $Z(U(u(2)_\h))$. Such a map $\al $ can be constructed
via the Harish-Chandra isomorphism. We introduce it on a larger algebra and by using another method.

Namely, we consider the following algebra
$$
\A=(\K(t,r)\otimes \Sym(su(2))) /\langle x^2+y^2+z^2-r^2\rangle,
$$
where $r$ is a new central generator and $\K(t,r)$ is the algebra of rational functions in $t$ and $r$. Thus, $r$ has the
meaning of the length of the radius vector of a point $(x,y,z)$.

Now, introduce the "quantum radius" by setting
 \be
\hat r={\hat \mu \over 2i}={\sqrt{\rule{0pt}{3.5mm}\h^2-4\Cas}\over 2i},
\label{rh}
\ee
where $\hat \mu $ is a new notation for the element introduced in (\ref{cas-mu}). Note that $\hat r$ coincides with the usual
radius provided ${\h}=0$.

This motivate the following quantum analog of the algebra $\A $
$$
\Ah=(\K(\hat t,\hat r) \ot U(su(2)_\h))/\langle \hat x^2+\hat y^2+\hat z^2+h^2-\hat r^2 \rangle
$$
where we put $\h=2ih$. Besides, if we assume $h$ to be real we get that $\h $ is purely imaginary and
the generators $\hat x, \hat y, \hat z$ can be represented as self-adjoint operators in a Hilbert space or a Verma module.

Now, define the quantizing map $\al:\A\to \Ah $ by the following relation
\be
\al(f(t,r) u)= f(\hat t,\hat r) \hat u,
\label{alpha-map}
\ee
where $\hat u$ is any element of the module $V^k$ above (see (\ref{sum})) and $u$ is its classical counterpart. Considering
the complexification of these modules and taking into account that $\al $ is $SU(2)$-invariant, it suffices to put
$$
u=b^k,\quad\hat u=\hat b^k,\quad {\rm where}\quad b = -(y+ix)\in\Sym(sl(2)),\quad\hat b = - (\hat y+i\hat x)\in U(sl(2)_\h).
$$

Now, extend the quantizing map on differential operators. With any given differential operators $\D=\sum_\beta a_\beta
\pa_\beta$ where $\beta$ is a multi-index, the sum is finite, and $a_\beta\in \A $, we associate the operator
\be
\al(\D):=\sum_\beta \al(a_\beta) \hat\pa_\beta
\label{alpha-oper}
\ee
with coefficients from $\Ah$.

\begin{definition}{\rm
The map $\D\mapsto \al(\D)$ is called the quantization of differential operators. }
\end{definition}

The operator $\al(\D)$ is acting on the algebra $U(u(2)_\h)$. However, we want it to be defined on the algebra $\Ah $.
The results of the previous section enable us to define such an action for the invariant operators (\ref{first}) and (\ref{second})
and consequently for any their combination with coefficients from $\Ah $. (Note that defining an action of the partial
derivatives on the whole algebra $\Ah $ is a more subtle deal and we leave it for subsequent publications.)

Similarly to formula (\ref{sum}) we present the algebra $\Ah $ as follows
$$
\Ah\cong \bigoplus_{k=0}^\infty\Bigl( \K(\hat t,\hat r)\ot V^k \Bigr).
$$

The components $\K(\hat t,\hat r)\ot V^k$ are also called {\em isotypic}.

According to proposition \ref{prop:20} each of the operators (\ref{first}) and (\ref{second}) maps such an isotypic component to itself.
Consequently, the same claim is valid for any linear combination of these operators with coefficients from $\Ah $. Below, we consider
an example of such an operator, namely, the Laplace-Beltrami operator corresponding to some rotationally symmetric metrics.   However, first, we construct the quantization of some invariant operators.

Consider the quantization of the operators $Q=x\pa_x+y\pa_y+z\pa_z$ and $Q^2$. First, we present the operator $Q^2$ in the ordered
form, moving all the partial derivatives to the right hand side position:
$$
Q^2 = Q + 2(xy\,\pa_x\pa_y +yz\,\pa_y\pa_z+zx\,\pa_z\pa_x) +x^2\pa^2_x+y^2\pa^2_y+z^2\pa^2_z.
$$
Then, we have to find the quantum images of the coefficients at the derivatives (see (\ref{alpha-oper})).
The definition (\ref{alpha-map}) of the quantizing map $\alpha $ gives
$$
\alpha(x) = \hat x,\qquad \alpha(x^2) = \hat x^2 +\frac{\h^2}{12},\qquad \alpha(xy) = \frac{\hat x\hat y + \hat y\hat x}{2}.
$$
The answer for other coefficients of $Q^2$ can be obtained from the above formulae by the cyclic substitution
$x\rightarrow y\rightarrow z $. Finally, we get the answer
$$
\hat Q =:\al(Q)= \hat x\hat \pa_x+\hat y\hat\pa_y+\hat z\hat\pa_z,
$$
$$
\al(Q^2)=\hat Q^2 + \frac{\h^2}{12}\,\hat \Delta -\frac{\h}{2} \hat Q \hat\pa_t.
$$
where $\hat \Delta = \hat\pa^2_x + \hat\pa^2_y+ \hat\pa^2_z$.
Emphasize that $\al(Q^2)\not= \hat Q^2$.

Now, we consider the quantization of a free massless scalar field in the space equipped with a
Schwarzschild-type metric
$$
\varphi(r) \, dt^2-\varphi(r)^{-1}\, dr^2-r^2 d\Omega^2.
$$
Here $d\Omega^2$ is the area form of the unit sphere and $\varphi(r)$ is a rational function. It is just Schwarzschild
metric provided $\varphi(r)=1-{r_g\over r}$.

The corresponding Laplace-Beltrami (LB) operator describing the dynamics reads
$$
\square_{LB}=\varphi(r)^{-1}\pa_t^2-\varphi(r)\pa_r^2- \frac{1}{r^2}\Bigl(X^2+Y^2+Z^2\Bigr)-{1\over r^2}
\pa_r(\varphi(r)\, r^2)\pa_r,
$$
where $X=y\, \pa_z-z\,\pa_y$, $Y=z\, \pa_x-x\,\pa_z$, $Z=x\, \pa_y-y\,\pa_x$. By using the relation
$$
\frac{1}{r^2}\Bigl(X^2+Y^2+Z^2\Bigr) =\De-\pa_r^2-{2\over r}\pa_r,
$$
we rewrite $\square_{LB}$ in the form
\be
\square_{LB}= \varphi(r)^{-1}\pa_t^2-(\varphi(r)-1)\pa_r^2 -\Bigl({2\over r} (\varphi(r)-1) +
\pa_r \varphi(r)\Bigr)\pa_r-\De.
\label{LB}
\ee

In order to quantize this operator we have to find quantum analogs of the operators $\pa_r$ and $\pa_r^2$.
Note, that in the classical setting
$$
\pa_r={1\over r}(x\,\pa_x+y\,\pa_y+z\,\pa_z)=\frac{1}{r}\,Q.
$$
As for the operator $\pa_r^2$, we have
$$
\pa_r^2={1\over r^2}\, (x^2\pa_x^2+y^2\pa_y^2+z^2\pa_z^2+2xy\pa_x\pa_y+2yz\pa_y\pa_z+2zx\pa_z\pa_x)=
\frac{1}{r^2}\,(Q^2-Q).
$$

Finally, formula (\ref{LB}) can be presented as
\be
\square_{LB}= \varphi(r)^{-1}\pa_t^2-\frac{\varphi(r)-1}{r^2}\,Q^2-\frac{1}{r}\,\Bigl(\frac{\varphi(r)-1}{r}+
\pa_r \varphi(r)\Bigr)Q -\De.
\label{LB1}
\ee
So, according to our scheme, the quantum counterpart $\hat \square_{LB} = \alpha(\square_{LB})$ of the operator $\square_{LB}$ is
$$
\hat\square_{LB}= \varphi(\hat r)^{-1}\hat\pa_t^2-\frac{\varphi(\hat r)-1}{\hat r^2}\,\hat Q^2-
\Bigl(\frac{\varphi(\hat r)-1} {\hat r^2}+\frac{(\pa_r\varphi)(\hat r)}{\hat r}\Bigr)\hat Q+
\frac{\h}{2} \frac{\varphi(\hat r)-1} {\hat r^2}\hat Q \hat\pa_t -
\Bigl(1+\frac{\h^2}{12}\,\frac{\varphi(\hat r)-1}{\hat r^2}\Bigr)\hat\De.
$$

Thus, the operator $\hat\square_{LB}$ is a linear combination of the operators (\ref{first}) and (\ref{second}) with coefficients
from $\K(t,\mu)$. Consequently, it is well defined on the whole algebra $\Ah $ and maps each isotypic component to itself.
More precisely, the result of the action $\hat\square_{LB}(f(\hat t, \hat r)\hat u)$, where $f$ is a rational function in two variables and $\hat u\in V^k$, is
an element of the form $g_k (\hat t, \hat r) \hat u$, where $g_k (\hat t, \hat r)$ is also a rational function.
Thus, the action of the operator $\hat\square_{LB}$
can be described by a series of the maps
$$
\K(\hat t, \hat r)\to \K(\hat t, \hat r):\qquad f(\hat t, \hat r)\mapsto g_k(\hat t, \hat r),\quad \forall\,k\in \NN.
$$

In conclusion, we want to emphasize that the variables $\hat t$ and $\hat r $, coming in the final version of the above
model, are commutative. So, this model can be treated by means of the commutative algebra and analysis. Thus, it would be
interesting to find the spectrum of the operator $\hat\square_{LB}$. However, in contrast with the classical case our
model is based on difference (not differential) operators. Also, note that our method is hopefully valid on the RE and
modified RE algebras. The only problem in this case consists in a reasonable introducing a quantum radius and extending
the action of invariant operators to rational functions. We plan to exhibit this construction in our subsequent publications.

\end{document}